\documentclass[11pt]{amsart}

\begin{document}
\newtheorem{theorem}{Theorem}[section]
\newtheorem{proposition}[theorem]{Proposition}
\newtheorem{lemma}[theorem]{Lemma}
\newtheorem{corollary}[theorem]{Corollary}
\newtheorem{defn}[theorem]{Definition}
\newtheorem{conjecture}[theorem]{Conjecture}

\theoremstyle{definition}
\newtheorem{definition}{Definition}

%\numberwithin{equation}{section}

\newcommand{\R}{\mathbb R}
\newcommand{\TT}{\mathbb T}
\newcommand{\Z}{\mathbb Z}
\newcommand{\HH}{\mathbb H}
\newcommand{\vE}{\mathcal E} % lattice points on the sphere
\newcommand{\dist}{\operatorname{dist}}

\title[Restriction of toral eigenfunctions to hypersurfaces]
{Restriction of toral eigenfunctions to hypersurfaces }
\author{Jean Bourgain and Ze\'ev Rudnick}
\address{School of Mathematics, Institute for Advanced Study,
Princeton, NJ 08540 }
\email{bourgain@ias.edu}
\address{Raymond and Beverly Sackler School of Mathematical Sciences,
Tel Aviv University, Tel Aviv 69978, Israel and
School of Mathematics, Institute for Advanced Study,
Princeton, NJ 08540 }
\email{rudnick@post.tau.ac.il}
\date{\today}
\begin{abstract}
Let $\TT^d = \R^d/\Z^d$ be the $d$-dimensional flat torus.
We establish for $d=2,3$ uniform upper and lower bounds on the restrictions of the eigenfunctions of
the Laplacian to smooth hyper-surfaces with non-vanishing curvature.

\end{abstract}

\maketitle

\section{Introduction and statements}

%This note has its origin in the work of Burq, G\'erard and Tzvetkov
%\cite{BGT} which establishes bounds for the $L^2$-norm of the
%restriction of the eigenfunctions of a smooth Riemannian surface $M$
%without boundary to smooth curves.

Let $M$ be a smooth Riemannian surface without boundary, $\Delta$
the corresponding Laplace-Beltrami operator and  $\Sigma$ a smooth
curve in $M$. Burq, G\'erard and Tzvetkov \cite{BGT} established
bounds for the $L^2$-norm of the restriction of eigenfunctions of
$\Delta$ to the curve $\Sigma$, showing that if $-\Delta
\varphi_\lambda = \lambda^2 \varphi_\lambda$, $\lambda>0$, then
\begin{equation}\label{BGT all}
 ||\varphi_\lambda||_{L^2(\Sigma)} \ll \lambda^{1/4} ||\varphi_\lambda||_{L^2(M)}
\end{equation}
and if $\Sigma$ has non-vanishing geodesic curvature then
\eqref{BGT all} may be  improved to
\begin{equation}\label{BGT curved}
 ||\varphi_\lambda||_{L^2(\Sigma)} \ll \lambda^{1/6} ||\varphi_\lambda||_{L^2(M)}
\end{equation}
%As pointed out in \cite{BGT},
Both \eqref{BGT all}, \eqref{BGT curved} are saturated for the sphere
$S^2$.

% one should expect substantially better bounds for hyperbolic surfaces.
%\begin{verbatim}
%One expects the supremum to be \lambda^\epsilon
%
%In fact one expects there to be a limit V for the L^2 norm
%of the restriction (rather than a bound)
%and that the value distribution of the eigenfunction restricted
%to the curve be Gaussian with mean zero and variance V.
%\end{verbatim}

In \cite{BGT} it is observed that for the flat torus $M=\TT^2$,
\eqref{BGT all} can be improved to
\begin{equation}\label{torus epsilon}
 ||\varphi_\lambda||_{L^2(\Sigma)} \ll \lambda^{\epsilon}
 ||\varphi_\lambda||_{L^2(M)},\quad \forall \epsilon>0
\end{equation}
due to the fact that there is a corresponding bound on the supremum
of the eigenfunctions.
They raise the question whether in \eqref{torus epsilon}  the factor
$\lambda^{\epsilon} $  can be replaced by a constant, that is whether
there is a uniform $L^2$ restriction bound.
As pointed out by Sarnak \cite{Sarnak}, if we take $\Sigma$ to be a
geodesic segment on the torus, %(that is a straight line segment),
this particular problem is essentially equivalent to the currently open question
of whether on the circle $|x|=\lambda$, the number of
lattice points  on an arc of size $ \lambda^{1/2}$ admits a uniform bound.
%Jarnik \cite{Jarnik} proved that an arc of length
%$\lambda^{1/3}$ contains at most two lattice points,
%and Cilleruelo and Cordoba \cite{CC} showed that for any $\delta<\frac 12$,
%arcs of length $\lambda^{\delta}$ contain at most $M(\delta)$ lattice points
%(in \cite{CG} it is conjectured that this remains true for any $\delta<1$).

In \cite{BGT} results similar to \eqref{BGT all} are also established
in the higher dimensional case for restrictions of eigenfunctions to smooth
submanifolds, in particular \eqref{BGT all} holds
for codimension-one submanifolds (hypersurfaces) and is sharp for the
sphere $S^{d-1}$.
Moreover, \eqref{BGT curved} remains valid for hypersurfaces
with nonvanishing curvature \cite{Hu}.

%in particular that \eqref{BGT curved} also holds for curved codimension-one hypersurfaces.

In this note we pursue the improvements of \eqref{BGT curved} for
the standard flat $d$-dimensional tori $\TT^d=\R^d/\Z^d$, considering the restriction to
(codimension-one) hypersurfaces $\Sigma$ with non-vanishing
curvature.
\begin{theorem}\label{thm 1}
 Let $d=2,3$ and let $\Sigma\subset \TT^d$ be a smooth hypersurface
 with non-zero curvature.
There are constants $0<c  <C  <\infty$ and $\Lambda >0$, all depending
on $\Sigma$, so that all eigenfunctions $\varphi_\lambda$ of the
Laplacian on $\TT^d$ with  $\lambda>\Lambda$ satisfy
\begin{equation}
 c ||\varphi_\lambda||_2 \leq ||\varphi_\lambda||_{L^2(\Sigma)} \leq C
 ||\varphi_\lambda||_2
\end{equation}
\end{theorem}

Observe that for the lower bound, the curvature assumption is necessary,
since the eigenfunctions $\varphi(x) = \sin(2\pi n_1 x_1)  $ all
vanish on the hypersurface $x_1=0$.
In fact this lower bound implies that a curved hypersurface cannot be
contained in the nodal set of eigenfunctions with arbitrarily large
eigenvalues.

The proof of Theorem~\ref{thm 1} (which will be sketched in the next section for the easy case of $d=2$)  permits also to introduce a notion of ``relative quantum limit'' for restrictions to $\Sigma$ as above, but we will not discuss this further here.

It is reasonable to believe that Theorem~\ref{thm 1} holds in any
dimension, and one could further conjecture an upper bound without
curvature assumptions. At this point, we may only state an improvement
of the exponent $1/6$: % in \eqref{BGT curved}:
\begin{theorem}\label{thm 2}
 For all $d\geq 4$ there is $\rho(d)<\frac 16$ so that if
 $\varphi_\lambda$ is an eigenfunction of the Laplacian on $\TT^d$,
%$-\Delta \varphi_\lambda = \lambda^2 \varphi_\lambda$,
and $\Sigma \subset \TT^d$ is a smooth compact hypersurface with positive curvature, then
\begin{equation}
  ||\varphi_\lambda||_{L^2(\Sigma)} \ll \lambda^{\rho(d)} ||\varphi||_2
\end{equation}
\end{theorem}

\section{Proof of Theorem~\ref{thm 1} for $d=2$}

Denote by $\sigma$ the normalized arc-length measure on the curve
$\Sigma$. Using the method of stationary phase, one sees that if
$\Sigma$ has non-vanishing curvature then the Fourier transform
$\widehat\sigma$ decays as
\begin{equation}\label{stationary phase}
 |\widehat\sigma(\xi) | \ll |\xi|^{-1/2}, \quad \xi \neq 0
\end{equation}
Moreover   $|\widehat\sigma(\xi)|\leq \widehat\sigma(0)=1$ with equality only for
$\xi=0$, hence
\begin{equation}\label{upper bd for sigma}
 \sup_{0\neq \xi\in \Z^2} |\widehat\sigma(\xi)| \leq 1-\delta
\end{equation}
for some $\delta=\delta_\Sigma>0$.

An eigenfunction of the Laplacian on $\TT^2$ is a trigonometric
polynomial of the form
\begin{equation}\label{fourier expansion}
 \varphi(x)  = \sum_{|n|=\lambda} \widehat\varphi(n)e(n\cdot x)
\end{equation}
(where $e(z):=e^{2\pi i z}$), all of whose frequencies lie in the set $\vE:=  \Z^2 \cap \lambda S^1  $.
As is well known, in dimension $d=2$, $\#\vE\ll \lambda^\epsilon$ for
all $\epsilon>0$. Moreover, by a result of Jarnik \cite{Jarnik}, any
arc on $\lambda S^1$ of length at most $c\lambda^{1/3}$ contains
at most two lattice points
%Jarnik \cite{Jarnik} proved that an arc of length
%$\lambda^{1/3}$ contains at most two lattice points,
(Cilleruelo and Cordoba \cite{CC} showed that for any $\delta<\frac 12$,
arcs of length $\lambda^{\delta}$ contain at most $M(\delta)$ lattice points
and in \cite{CG} it is conjectured that this remains true for any $\delta<1$).
Hence we may partition
\begin{equation}
 \vE = \coprod_\alpha \vE_\alpha
\end{equation}
where $\#\vE_\alpha \leq 2$ and $\dist(\vE_\alpha, \vE_\beta)>c
\lambda^{1/3}$ for $\alpha\neq \beta$.
Correspondingly we may write
\begin{equation}
 \varphi = \sum_\alpha \varphi^\alpha,\quad
\varphi^\alpha(x)  =\sum_{n\in \vE_\alpha} \widehat\varphi(n)e(nx)
\end{equation}
so that $||\varphi||_2^2 = \sum_\alpha ||\varphi^\alpha||_2^2$  and
\begin{equation}
 \int_\Sigma |\varphi|^2 d\sigma = \sum_\alpha \sum_\beta \int_\Sigma
 \varphi^\alpha \overline{\varphi^\beta} d\sigma
\end{equation}

Applying \eqref{stationary phase}  we see that  $\int_\Sigma
\varphi^\alpha \overline{\varphi^\beta} d\sigma \ll \lambda^{-1/6}$ if
$\alpha\neq \beta$ and because $\#\vE \ll \lambda^\epsilon$ the total
sum of these nondiagonal terms is bounded by
$\lambda^{-1/6 +\epsilon} ||\varphi||_2^2$.
It suffices then to show that the diagonal terms satisfy
\begin{equation}\label{diagonal term}
 \delta ||\phi^\alpha||_2^2 \leq \int_\Sigma |\phi^\alpha|^2 d\sigma
 \leq 2 ||\phi^\alpha||_2^2
\end{equation}
This is clear if $\vE_\alpha=\{n\}$ while if $\vE_\alpha=\{m,n\}$ then
\begin{equation}
\int_\Sigma |\phi^\alpha|^2 d\sigma   =
|\widehat\varphi(m)|^2 + |\widehat\varphi(n)|^2 +
2\mbox{Re }\widehat\varphi(m)\overline{\widehat\varphi(n)} \widehat\sigma(m-n)
\end{equation}
and then \eqref{diagonal term} follows from \eqref{upper bd for sigma}.
Thus we get  Theorem~\ref{thm 1} for  $d=2$.

\section{The higher-dimensional case}
The proof of Theorem~\ref{thm 1} for dimension $d=3$   is
considerably more involved. Arguing along the lines of the
two-dimensional case gives an upper bound of $\lambda^{\epsilon}$.
To get the uniform bound of Theorem~\ref{thm 1} for $d=3$ and the
results of Theorem~\ref{thm 2}, we need to replace the upper bound
\eqref{stationary phase} for the Fourier transform of the
hypersurface measure by an asymptotic expansion, and then exploit
cancellation in the resulting exponential sums over the sphere. A
key ingredient there is controlling the number of lattice points in
spherical caps.

\subsection{Distribution of lattice points on spheres}
To state some relevant results, denote as before by $\vE=\Z^d\cap \lambda S^{d-1}$
the set of lattice points on the sphere of radius $\lambda$. %(where we assume $\lambda^2$ is an integer).
We have $\#\vE \ll \lambda^{d-2+\epsilon}$.
Let $F_d(\lambda, r)$ be the maximal number of lattice points in the intersection of $\vE$
with a spherical cap of size $r>1$.
A higher-dimensional analogue of Jarnik's theorem %\cite{Andrews}
implies that if $r\ll \lambda^{1/(d+1)}$
then all lattice points in such  a cap are co-planar, hence $F_d(r,\lambda)\ll r^{d-3+\epsilon}$
in that case, for any $\epsilon>0$. For larger caps, we show:
\begin{proposition}
i) Let $d=3$. Then for any $\eta<\frac 1{15}$,
 \begin{equation} \label{bd for F3}
  F_3(\lambda,r) \ll \lambda^\epsilon \left(   r (\frac r\lambda)^\eta +1 \right)
 \end{equation}

ii) Let $d=4$. Then
\begin{equation}
 F_4(\lambda,r) \ll \lambda^\epsilon \left( \frac{r^3}\lambda + r^{3/2} \right)
\end{equation}

iii) For $d\geq 5$ we have
\begin{equation}
 F_d(\lambda,r) \ll \lambda^\epsilon \left( \frac{r^{d-1}}\lambda + r^{d-3} \right)
\end{equation}
(the factor $\lambda^\epsilon$ is redundant for large $d$).
\end{proposition}
The term $r^{d-1}/\lambda$ concerns the equidistribution of $\vE$,
while the term $r^{d-3}$ measures deviations related to accumulation in lower dimensional strata.

The second result expresses a mean-equidistribution property of $\vE$.
Partition the sphere $\lambda S^{2}$ into sets $C_\alpha$ of size $\lambda^{ 1/2}$,
for instance by intersecting with cubes of that size. Since $\#\vE\ll \lambda^{1+\epsilon}$, one may expect that
$\#C_\alpha\cap \vE\ll \lambda^\epsilon$.
Using Siegel's mass formula for the number of representations of an integral quadratic form by the genus of a quadratic form,
we show (in joint work with P. Sarnak \cite{BRS}) that this holds in the mean square:
\begin{proposition}\label{prop:siegel}
 %For $d\geq 3$ we have
\begin{equation}
 \sum_\alpha [\#(\vE\cap C_\alpha)]^{2} \ll \lambda^{1+\epsilon}, \quad \forall \epsilon>0
\end{equation}
\end{proposition}

\subsection{Exponential sums on the sphere}

Let $1<r<\lambda$ and let $C$, $C'$ be spherical $r$-caps on $\lambda S^{d-1}$ of mutual distance at least $10r$. Following the argument for $d=2$, we need to bound exponential sums of the form
\begin{equation}\label{double sum}
 \sum_{n\in C} \sum_{n'\in C'} \widehat \varphi(n) \overline{\widehat\varphi(n')} e(h(n-n')) , \qquad ||\varphi||_2=1
\end{equation}
where $h$ is the support function of the hyper-surface $\Sigma$, which appears in the asymptotic expansion of the Fourier transform of the surface measure on $\Sigma$, see \cite{Herz}. For instance, in the case that $\Sigma=\{|x|=1\}$ is the unit sphere then $h(\xi) = |\xi|$.

Consider from now on the case  $d=3$. For $r<\lambda^{1-\epsilon}$
we simply estimate \eqref{double sum} by $F_3(\lambda,r)$ (see
\eqref{bd for F3}). When $\lambda^{1-\epsilon} <r<\lambda$ this
bound does not suffices and we need to exploit cancellation in the
sum \eqref{double sum}.

\begin{proposition}
 There is $\delta>0$ so that \eqref{double sum} admits a bound of $\lambda^{1-\delta}$ for $\lambda\gg 1$.
\end{proposition}
This statement depends essentially on the equidistribution of $\vE$ in $\sqrt{\lambda}$-caps,
as expressed in Proposition~\ref{prop:siegel}. % for $d=3$.

We finally formulate an example of a bilinear estimate involved in analyzing \eqref{double sum}.
\begin{proposition}
 Let $\beta\gg 1$ and $X,Y\subset [0,1]$ arbitrary discrete sets such that $|x-x'|, |y-y'|>\beta^{-1/2}$
for $x\neq x'\in X$ and $y\neq y'\in Y$. Then
\begin{equation}
 \left| \sum_{x\in X}\sum_{y\in Y} e( \beta xy + \beta^{1/3}x^2y^2) \right| \ll \beta^{23/24+\epsilon}
\end{equation}
for all $\epsilon>0$.
\end{proposition}
Note that the nonlinear term in the phase function is crucial for a non-trivial bound to hold in this setting.

\subsection{Acknowledgements and grant support} We thank Peter Sarnak for several stimulating discussions.
Supported in part by N.S.F. grant DMS 0808042 (J.B.)  and Israel Science Foundation  grant No. 925/06 (Z.R.).

\end{document}